\documentclass[11pt]{amsart}

\usepackage{amsmath,amsfonts,amsthm,amscd,amssymb,graphicx}
\usepackage[utf8]{inputenc}
\usepackage{graphicx}

\usepackage[text={6.5in,9in},centering,includefoot,foot=0.6in]{geometry}

\usepackage{latexsym,verbatim}
\usepackage{epsfig,epstopdf,color}

\newtheorem{Lemma}{Lemma}[section]

\newtheorem{Proposition}[Lemma]{Proposition}
\newtheorem{Remark}[Lemma]{Remark}
\newtheorem{Theorem}{Theorem}

\newenvironment{Proof}[1][.]%
 {\begin{trivlist}\item[]\textbf{Proof#1 }}%
 {\hspace*{\fill}$\rule{0.3\baselineskip}{0.35\baselineskip}$\end{trivlist}}

 {\begin{trivlist}\item[]\textbf{Acknowledgments }}{\end{trivlist}}

\makeatletter\@addtoreset{equation}{section}\makeatother

\def\Re{\mathop\mathrm{Re}\nolimits}    

\newcommand{\rmd}{\mathrm{d}}           
\newcommand{\rmi}{\mathrm{i}}           

\newcommand{\uu}{{\bf u}}
\newcommand{\vv}{{\bf v}}

\newcommand{\cLa}{{ {\mathcal L}^a }}
\newcommand{\cLell}{{ {\mathcal L}_{\ell} }}
\newcommand{\cLaell}{{ {\mathcal L}^a_{\ell} }}
\newcommand{\hwell}{{ \hat{w}_{\ell} }}
\newcommand{\huell}{{ \hat{u}_{\ell} }}

\newcommand{\hw}{{ \hat{w} }}
\newcommand{\Bell}{{ B^{\ell} }}
\newcommand{\Cell}{{ C^{\ell} }}

\begin{document}

\title[Metastability and Navier-Stokes]{Metastability and rapid convergence to quasi-stationary bar states for the 2D Navier-Stokes Equations}

\author{ Margaret Beck }
\address{Department of Mathematics and Statistics, Boston University, Boston, USA}
\curraddr{Department of Mathematics\\ Heriot-Watt University\\  Edinburgh, Scotland}
\email{M.Beck@hw.ac.uk}
\author{C. Eugene Wayne}
\address{Department of Mathematics and Statistics, Boston University, Boston, USA}
\email{cew@bu.edu}

\date{\today}

\maketitle

\begin{abstract}
Quasi-stationary, or metastable, states play an important role in two-dimensional turbulent fluid flows 
where they often emerge on time-scales
much shorter than the viscous time scale, and then dominate the dynamics
for very long time intervals.  In this paper we propose a dynamical systems
explanation of the metastability of an explicit family of solutions, referred to as bar states, of the two-dimensional 
incompressible Navier-Stokes equation on the torus.  These states are physically relevant because 
they are associated with certain maximum entropy solutions of the Euler equations, and they have been 
observed as one type of metastable state in numerical studies of two-dimensional turbulence. For small 
viscosity (high Reynolds number), these states are quasi-stationary in the sense that they decay on the 
slow, viscous timescale. Linearization about these states leads to a time-dependent operator.
We show that if we approximate this operator by dropping a higher-order, non-local term, it produces
a decay rate much faster than the viscous decay rate.  We also provide numerical evidence that 
the same result holds for the full linear operator, and that our theoretical results give the optimal decay
rate in this setting.  
\end{abstract}


\section{Introduction}

It has been observed, both numerically and experimentally, that solutions to the nearly-inviscid two-dimensional incompressible Navier-Stokes equation on the torus will rapidly approach certain long-lived quasi-stationary states \cite{Couder83, MatthaeusStriblingMartinez91, McWilliams84, McWilliams90}. Evidence also suggests that these quasi-stationary, or metastable, states are connected with stationary solutions of the inviscid Euler equations. In \cite{YinMontgomeryClercx03}, the authors argue that certain maximum entropy solutions, determined using numerical and analytical techniques, of the inviscid Euler equation are the most probable quasi-stationary states that one would observe. They refer to these states as dipole and bar states and confirm their predictions by numerically analyzing the Navier-Stokes equation and demonstrating that most solutions do indeed converge rapidly to one of these two classes of states. Furthermore, in \cite{BouchetSimonnet09}, the authors consider a stochastically forced Navier-Stokes equation, and observe both the dipole and the bar states as statistical equilibria. 

There are explicit exact solutions of the Navier-Stokes equation that qualitatively match the bar and dipole states. (See equations (\ref{E:bar}) and (\ref{E:dipole}), below.) The bar states are often referred to as Kolmogorov flow \cite{BouchetMorita10}, and the dipole as a Taylor-Green vortex  \cite[pg. 172]{MajdaBertozzi02}, but we adopt the terminology of  \cite{YinMontgomeryClercx03} when referring to these exact solutions. In this work, we provide analytical and numerical results indicating that there is a class of initial data for which solutions to the linearized equation rapidly converge to the bar states.   More precisely, for small viscosity $0 < \nu \ll 1$,
we prove that for a natural approximation to the linearized operator,  there exists an invariant subspace for the linearized equation in which solutions converge to the bar states at a rate that is $\mathcal{O}(e^{-\sqrt{\nu}t})$, which is much faster than the decay of the bar states themselves, which is $\mathcal{O}(e^{-\nu t})$. Furthermore, we present the results of numerical computations
that indicate that this decay rate is optimal, not just for our approximation but also for the complete linearization
of the Navier-Stokes equation about the bar states.
We do not analyze the dipole solutions here - see \S \ref{S:discussion} for more information about these states. 

This work can be viewed as an extension of our prior research on metastable behavior in Burgers equation, \cite{BeckWayne11}.  In that work we showed that corresponding to each long-time asymptotic state of Burgers equation there is a one-dimensional, invariant manifold of viscous $N$-waves which attracts nearby solutions at a rate much faster than the viscous time scale and which governs the behavior of solutions for very long times, before they reach their
ultimate, asymptotic state.  Furthermore, this manifold is tangent, at the fixed point representing the asymptotic state, to the eigenspace of the eigenvalue of smallest real part.
For the two-dimensional Navier-Stokes equation on the torus, the asymptotic state is the rest
solution in which the velocity of the fluid is zero.  Linearizing about this state, the eigendirections
corresponding to the smallest non-zero eigenvalues correspond
to the bar states, and these invariant manifolds (in this case, actually invariant subspaces) form attractive, invariant manifolds of metastable states, just as in the case of Burgers equation.

It is interesting to note, however, that the mechanism causing the separation in times scales for Burgers equation and the bar states of the Navier-Stokes equation appears to be different. Linearization about the bar states produces small eigenvalues that scale with
that square root of  the viscosity parameter. In contrast, linearization about the asymptotic states for Burgers equation leads to an operator whose eigenvalues are independent of viscosity, and the separation in time scales appears to result from exponentially large coefficients in the eigenfunction expansion. Both phenomena, however, are connected with the fact that neither linearized operator is self-adjoint.

Consider the two-dimensional incompressible Navier-Stokes equation on the torus, 
\[
\partial_t \uu = \nu \Delta \uu - \uu \cdot \nabla \uu + \nabla p, \qquad \nabla \cdot \uu = 0,
\]
where $-\pi < x < \pi$,  $-\pi < y < \pi$, $\uu: \mathbb{T}^2 \to \mathbb{R}^2$, $\uu(-\pi,y,t) = \uu(\pi, y ,t)$ and $\uu(x, -\pi, t) = \uu(x, \pi,t)$ for all $t \geq 0$. The pressure $p(x,y,t)$ is a scalar valued function that also satisfies the periodic boundary conditions, and the viscosity $\nu$ is assumed to be small: $0 < \nu \ll 1$. The fluid vorticity $\omega$ is defined in terms of the velocity $\uu$ via $\omega = \nabla \times \uu = \partial_x \uu_2 - \partial_y \uu_1$. Taking the curl of the above equation, we find the two-dimensional vorticity equation
\begin{equation} \label{eq:2Dvort}
\partial_t \omega = \nu \Delta \omega - \uu \cdot \nabla \omega,
\end{equation}
where the scalar-valued function $\omega$ also satisfies periodic boundary conditions. We assume without loss of generality that solutions to this equation have zero mean -- equation (\ref{eq:2Dvort}) conserves mass and, in fact, the periodic boundary conditions for $\uu$ imply that $\omega$ has zero mean. In this case, one can determine the velocity from the vorticity using the Biot-Savart law as follows. Given a periodic function $\omega \in L^2(\mathbb{T}^2)$, we can write it as its Fourier series
\[
\omega(x,y) = \sum_{k,l \in \mathbb{Z}} \hat{\omega}(k,l) e^{i(kx+l y)},
\]
where
\[
\hat{\omega}(k,l) = \frac{1}{4 \pi^2} \int_{[-\pi,\pi]\times[-\pi,\pi]} \omega(x,y) e^{-i(kx+l y)} dx dy.
\]
Using the fact that $\omega = \nabla \times \uu$ and the incompressibility condition $\nabla \cdot \uu = 0$, we have the Biot-Savart law
\[
\hat{\uu}(k,l) = \frac{i}{k^2+l^2}(l, -k) \hat{\omega}(k,l), \qquad (k,l) \neq (0,0).
\]
Alternatively, one can obtain this relationship using the fact that, if $\psi$ is the stream function, then $\Delta \psi = w$ and $\uu = \nabla^\perp \psi = ( - \partial_y \psi, \partial_x \psi)$. Since $w$ has mean zero, 
\[
\uu = (-\partial_y \Delta^{-1} w, \partial_x \Delta^{-1} w).
\]

The explicit family of solutions of (\ref{eq:2Dvort}) that qualitatively matches the maximum entropy solutions found in \cite[Figure 2c]{YinMontgomeryClercx03} and referred to as ``bar states" 
has vorticity given by
\begin{equation}\label{E:bar}
\omega^{b,m}(x,y,t) = e^{-\nu m^2 t} \cos(mx), \quad \tilde{\omega}^{b,m}(x,y,t) = e^{-\nu m^2 t} \sin(mx),
\end{equation}
and the associated velocity solutions are
\[
\uu^{b,m}(x,y,t) =  \frac{1}{m} e^{-\nu m^2 t} \begin{pmatrix} 0 \\ \sin mx \end{pmatrix}, \qquad \tilde{\uu}^{b,m}(x,y,t) = - \frac{1}{m} e^{-\nu m^2 t} \begin{pmatrix} 0 \\ \cos mx \end{pmatrix}.
\]
Note that the bar states can also be defined to be functions of $y$, rather than $x$, and any linear combination of $x$-bar states, or any linear combination of the $y$-bar states, is again an exact solution. 

We will focus primarily on these states when $m = 1$ and denote $\bar{\omega} = \omega^{b,1}$. When $0 < \nu \ll 1$, these states are quasi-stationary, because they evolve in time only on the very long timescale, $\mathcal{O}(e^{-\nu t})$. They are stationary states for the inviscid Euler equation (i.e. when $\nu = 0)$. Due to incompressibility and periodic boundary conditions, for any solution of (\ref{eq:2Dvort}) we have
\[
\frac{d}{dt}\frac{1}{2} \int_{\mathbb{T}^2} \omega^2(x,y) \rmd x\rmd y = -\nu \int_{\mathbb{T}^2} |\nabla \omega(x,y)|^2 \rmd x\rmd y \leq -C \nu \int_{\mathbb{T}^2} \omega^2(x,y) \rmd x\rmd y,
\]
where we have used Poincar\'e's inequality to obtain the final inequality. Thus, the slow decay rate of the bar states occurs on what is arguably the natural timescale for the system, determined by the viscosity $\nu$. 

Consider a solution of (\ref{eq:2Dvort}) of the form $\omega = a\bar{\omega} + w$, where the constant $a$ represents the amplitude of the bar state around which we will linearize. Since equation (\ref{eq:2Dvort}) conserves mass, we will assume that $\omega$, and hence also $w$, has mean zero. The perturbation then satisfies
\begin{eqnarray}
\partial_t w &=& \mathcal{L}(t)w - \uu \cdot \nabla w \label{E:lin_b} \\
\mathcal{L}(t) w &=& \nu \Delta w - ae^{-\nu t} [ \sin x \partial_y (1 + \Delta^{-1})] w, \nonumber
\end{eqnarray}
where $\uu$ is the velocity corresponding to $w$ via the Biot-Savart law. 

Our goal is to show that there exists an invariant subspace for the linear part
of equation (\ref{E:lin_b})  on which solutions decay very rapidly, at a 
rate that is $\mathcal{O}(e^{-\sqrt{\nu}t})$. See Theorem \ref{Thm:bar_conv}.
There are a couple of issues which have to be overcome to obtain this result. First, ignoring for the moment the time-dependence of the operator in (\ref{E:lin_b}), it must be shown that the addition of the term $\sin x \partial_y (1 + \Delta^{-1})$ to the Laplacian produces increased decay.  In the current paper we prove
this for an approximation to the linear operator which ignores the 
nonlocal part of this operator. This is done using the so-called hypocoercivity properties of the operator 
as developed in \cite{Villani09}. We also discuss the changes that will be needed to extend this result to the full linearized operator.  
We note that a similar approach has recently been used to study the 
linearization of the two-dimensional Navier-Stokes equation around an Oseen
vortex \cite{Deng12}.

Second, it's clear that we can't expect decay at the rate $\mathcal{O}(e^{-\sqrt{\nu}t})$ on the entire space. This is because there are explicit solutions of the linear part of equation (\ref{E:lin_b}) given by $e^{-\nu m^2  t \pm \rmi m x}$ and $e^{-\nu t \pm \rmi y}$. Therefore, we need to find a way to project off of these anomalous modes. If the operator was time-independent, this could be done using a spectral projection via the adjoint eigenfunctions associated with the eigenvalues $-\nu$ and $-m^2 \nu$. However, that will not work in this case. Nevertheless, we are able to find an invariant subspace in which solutions are rapidly decaying. We first
examine this subspace for the full linear operator, and then show that for our approximation, in which we omit the 
nonlocal term, the rapidly decaying modes have a simple expression in terms of the Fourier coefficients of the solution.

A further consequence of the time-dependence of the operator is that we cannot expect a result that holds for all time. It is the term $ae^{-\nu t} [ \sin x \partial_y (1 + \Delta^{-1})]$ that must be producing the rapid decay, and this term is $\mathcal{O}(1)$ only for times up to and including $\mathcal{O}(1/\nu)$. Hence, Theorem \ref{Thm:bar_conv} is only valid on a time interval of this size.

The reason that we focus on the bar states, rather than the dipoles, is primarily that these states seem easier to analyze from a mathematical point of view. Although in \cite{YinMontgomeryClercx03} it was found that a large class of initial conditions have the dipole as their quasi-stationary state, we note that there are certain situations in which bar states seem to be more physically relevant. For example, in \cite{BouchetSimonnet09} the authors find that the bar states become dominant if the torus is elongated slightly in either the $x$ or the $y$ direction. 

There has been much previous work related to the multiple time-scales in the evolution of solutions to the two-dimensional incompressible Navier-Stokes equation. Those most relevant for the present analysis have been referred to above. We do not attempt to provide a comprehensive list, but we mention a few other works that are in the same spirit, in the sense that they approach the problem from a dynamical systems point of view. In particular, in \cite{GallayWayne02,GallayWayne05} the authors analyze the global stability of the single Oseen vortex for the 2D Navier-Stokes equation on the plane using invariant manifolds and a Lyapunov functional, and in \cite{Gallay11_arma} the author provides rigorous asymptotic expansions for the slow motion of localized vortices, again for the 2D Navier-Stokes equation on the plane. The case of two identical vortices is considered in \cite{Gallay11}. The paper of Lan and Li \cite{LanLi08} considers eigenvalues of the linearization of the Navier-Stokes equation for small
values of the viscosity but does not examine the viscosity dependence in the way we do here, nor do they consider the relationship of this problem to metastable behavior.  We also note that there is potentially a connection between the $\sqrt{\nu}$ scaling, found here, and Taylor dispersion \cite{Taylor53}, although we have not explored this in detail. 

The rest of the paper is organized as follows. In \S \ref{S:anomolous_modes} we identify
a small set of anomalous modes that decay on the viscous timescale $\mathcal{O}(e^{-\nu t})$
 and define an invariant subspace containing none of these modes.
 In \S\ref{S:hypocoercivity}, we prove that for a slightly simplified approximation to this linear operator,
 solutions of the associated evolution equation decay at the rate of at least $\mathcal{O}(e^{-\sqrt{\nu} t})$. 
 We also discuss the possible extension of these results to the full linear equation.
 Next, in \S\ref{S:numerics} we present numerical evidence that the decay rate derived in the prior section is in fact optimal
 by studying the way in which the eigenvalues of the linearized operator depend on the viscosity, $\nu$.
 Finally, we conclude with a discussion and mention directions for future work in \S \ref{S:discussion}.


\section{Anomalous modes} \label{S:anomolous_modes}

Our main goal of this paper is to show that the bar states attract nearby solutions at a rate much faster than expected from the viscous
time scale.  However, there are some ``anomalous'' modes which do decay only at the viscous rate.  In the present section we
identify the set of states with this anomalous behavior and  show that it is ``small'' in the sense that if we
expand solutions in a Fourier series with respect to the basis $\{ e^{ikx} e^{ily}\}$, then no modes with $| l | >1$ behave
anomalously, and only a very special set of those with $| l | =1$ do.

The linear part of equation (\ref{E:lin_b}) is given by
\begin{equation} \label{E:linear}
w_t = \mathcal{L}(t)w, \quad \mathcal{L}(t) =  \nu \Delta w - ae^{-\nu t} [ \sin x \partial_y (1 + \Delta^{-1})] w.
\end{equation}
In Fourier space, this operator becomes
\begin{eqnarray}
\hat{\mathcal{L}}(t) \hat{w}(k,l) &=& - \nu(k^2+l^2) \hat{w}(k,l) -  \frac{l}{2}ae^{-\nu t} \Bigg[ \left( 1 - \frac{1}{(k-1)^2+l^2}\right) \hat{w}(k-1,l) \nonumber \\
&\qquad& \qquad - \left( 1 - \frac{1}{(k+1)^2+l^2}\right) \hat{w}(k+1,l) \Bigg].
\label{E:bar_lin}
\end{eqnarray}
We wish to show that solutions to equation (\ref{E:linear}) decay at a rate which is faster than $\mathcal{O}(e^{-\nu t})$. As mentioned above, there is one key obstacle preventing this from being true, in general, and that is the presence solutions of the form $e^{-\nu m^2 t}e^{\pm \rmi m x}$ and $e^{-\nu t}e^{\pm \rmi y}$, which we refer to as anomalous modes. Therefore, any solution whose Fourier series has $\hat{\omega}(\pm m, 0)(t) \neq 0$ and/or $\hat{\omega}(0, \pm1)(t) \neq 0$ cannot decay faster than $\mathcal{O}(e^{-\nu m^2 t})$. 

In some sense these solutions correspond to the eigenvalues $-\nu m^2$ for the linear operator. However, since the linear operator is time-dependent, these are not eigenvalues in the usual sense. Nevertheless, we will see they have a similar effect, and we will slightly abuse this terminology. If the operator were time-independent, we could simply determine the corresponding adjoint eigenfunctions and use them to construct spectral projections and an invariant subspace without any anomalous modes. However, in this time-dependent case we must do something slightly different.

\begin{Lemma} The functions $e^{-\nu m^2 t}e^{\pm \rmi m x}$, $m \in \mathbb{N}$, and $e^{-\nu t}e^{\pm \rmi y}$ are solutions of (\ref{E:linear}). In addition, $e^{-\nu m^2 t}e^{\pm \rmi m x}$ are solutions of the corresponding adjoint equation $\partial_t w = \mathcal{L}^*(t) w$, where $\mathcal{L}^*$ is the $L^2$ adjoint of $\mathcal{L}$, but $e^{-\nu t}e^{\pm \rmi y}$ are not.
\end{Lemma}

\begin{Proof}
This follows from a direct calculation, using the fact that
\[
\mathcal{L}^*(t) = \nu \Delta + ae^{-\nu t} [ (1 + \Delta^{-1}) \partial_y\sin x ].
\]
The key observation is that $\partial_y e^{\pm m \rmi x} = (1 + \Delta^{-1}) e^{\pm \rmi y} = 0$.
\end{Proof}

Because the functions $e^{-\nu m^2 t}e^{\pm \rmi m x}$ also solve the adjoint problem, the subspace 
\[
\mathcal{M}_{\mathrm{rapid}, x} := \{ \omega \in L^2(\mathbb{T}^2): \hat{\omega}(m, 0) = 0, \mbox{ } m \in \mathbb{Z} \}
\]
is invariant for the linear PDE (\ref{E:linear}). Therefore, we can avoid these anomalous modes by working only in this subspace. Unfortunately, it is less straightforward to remove the other anomalous modes, because the corresponding adjoint ``eigenfunctions" are not time-independent, as the following calculation shows.

Suppose we look for solutions to $-\nu \omega = \mathcal{L}^* \omega$, in attempt to determine a projection off the modes $e^{\pm \rmi y}$. In Fourier space, this equation is
\[
0 = \nu[1 - (k^2 + l^2)] \hat{\omega}(k,l) + \frac{l}{2}ae^{-\nu t}\left( 1 - \frac{1}{k^2+l^2}\right) [ \hat{w}(k-1,l) -  \hat{w}(k+1,l)].
\]
As mentioned above, $e^{\pm \rmi m x}$ are solutions. Are there others? When $(k,l) = (0, \pm 1)$, the above equation is always satisfied. Thus, without loss of generality we can take $\hat{\omega}(0, \pm 1) = 1$. However, we then find for $(k,\pm 1) = (1,\pm 1)$
\[
0 = -\nu \hat{\omega}(1,\pm 1) \pm \frac{1}{4}ae^{-\nu t}[ \hat{w}(0,\pm 1) -  \hat{w}(2,\pm 1)],
\] 
and for $(k,\pm 1) = (-1,\pm 1)$
\[
0 = -\nu \hat{\omega}(-1,\pm 1) \pm \frac{1}{4}ae^{-\nu t}[ \hat{w}(-2,\pm 1) -  \hat{w}(0,\pm 1)],
\] 
and so on. Thus, the effect of the $(0, \pm 1)$ mode will propagate out requiring nonzero, and time dependent, entries in other components of the adjoint eigenfunction. Formal calculations with this equation lead one to the following result. 

Define
\begin{equation} \label{E:def_pq}
p_j^\pm := \hat{w}(2j, \pm 1) + \hat{w}(-2j, \pm 1), \qquad q_j^\pm := \hat{w}(2j+1, \pm 1) - \hat{w}(-(2j+1), \pm 1),
\end{equation}
for $j \in \mathbb{N}$ and 
\begin{equation} \label{E:inv_anom}
\mathcal{M}_{\mathrm{rapid}, y} := \{ w \in L^2(\mathbb{T}^2): p_j^\pm = q_j^\pm = 0 \mbox{ for all } j \in \mathbb{N}\}.
\end{equation}

\begin{Proposition} \label{prop:anom} A solution $w(t) \in L^2(\mathbb{T}^2)$ of (\ref{E:linear}) satisfies $\hat{w}(0, \pm 1)(t) = 0$ for all $t \geq 0$ if and only if $w(0) \in \mathcal{M}_{\mathrm{rapid}, y}$, defined in (\ref{E:inv_anom}).
\end{Proposition}

\begin{Proof} Using equation (\ref{E:linear}), we find
\begin{eqnarray*}
\partial_t p_j^\pm &=& - \nu (4j^2 + 1)p_j^\pm \mp \frac{a}{2}e^{\nu t} \left[ \left( 1 - \frac{1}{(2j-1)^2 + 1}\right) q_{j-1}^\pm  - \left( 1 - \frac{1}{(2j+1)^2 + 1}\right) q_j^\pm \right] \\
\partial_t q_j^\pm &=& - \nu ((2j+1)^2 + 1)q_j^\pm \mp \frac{a}{2}e^{\nu t} \left[ \left( 1 - \frac{1}{(2j)^2 + 1}\right) p_j^\pm  - \left( 1 - \frac{1}{(2j+2)^2 + 1}\right) p_{j+1}^\pm \right].
\end{eqnarray*}
As a result, the space $\mathcal{M}_{\mathrm{rapid}, y}$ is invariant, so if $w(0) \in \mathcal{M}_{\mathrm{rapid}, y}$ then $\hat{w}(0, \pm 1)(t) = 0$ for all $t \geq 0$. 

To prove the other direction, suppose $\hat{w}(0, \pm 1)(t) = 0$ for all $t \geq 0$. Using the Fourier representation of $\mathcal{L}(t)$, given in (\ref{E:bar_lin}), we see that the dynamics of $\hat{w}(0, \pm 1)$ are determined by
\[
\partial_t \hat{w}(0, \pm 1) = -\nu \hat{w}(0, \pm 1) \mp \frac{1}{4}Ae^{-\nu t}[  \hat{w}(-1,\pm 1)  - \hat{w}(1,\pm 1)].
\]
Thus, in order to have $\hat{w}(0, \pm 1)(t) = \frac{1}{2} p_0^\pm = 0$ for all $t \geq 0$, it must be the case that $\hat{w}(-1,\pm 1)  - \hat{w}(1,\pm 1) = -q_0^\pm = 0$ for all $t \geq 0$. Looking now at the above equations for $p_j^\pm$ and $q_j^\pm$, we have
\[
\partial_t q_0^\pm = - 2\nu q_0^\pm \mp \frac{A}{2}e^{\nu t} \left[  - \frac{4}{5}p_1^\pm \right],
\]
which then implies we must have $p_1^\pm = 0$ for all time. One can continue in this manner to show that $q_j^\pm = p_j^\pm = 0$. Another way to think about this is that, if one were to define $u^\pm = (p_0^\pm, q_0^\pm, p_1^\pm, q_1^\pm, \dots)$, then $u_t^\pm = A(t) u^\pm$, where $A(t)$ is a tridiagonal matrix. The only solution to this linear equation with $u^\pm(t) = (0, 0, p_1^\pm(t), q_1^\pm(t), \dots)$ is the solution that is identically zero. This proves the result.
\end{Proof}

As a result, we define
\begin{equation} \label{E:defM}
\mathcal{M} := \mathcal{M}_{\mathrm{rapid}, x} \cap \mathcal{M}_{\mathrm{rapid}, y},
\end{equation}
which is the subspace one which one can hope to 
have  rapid decay of solutions.  In \cite{YinMontgomeryClercx03} the authors found rapid convergence to the bar states only for  specific types of initial data. (More general initial data converge rapidly to the dipole solutions -- see \S \ref{S:discussion}.) Examples of such initial conditions that they found do indeed seem to match the condition imposed by the subspace $\mathcal{M}$. See \cite[Figures 8 \& 13]{YinMontgomeryClercx03}. Note, however, that in spite of the restrictions the
anomalous modes impose, the bar states are still physically relevant due to the fact that they can be viewed as maximum entropy solutions of the Euler equation, as well as the fact that they seem to dominate the flow when the aspect ratio of the torus is changed slightly \cite{BouchetSimonnet09}. 


\section{Rapid decay via hypocoercivity} \label{S:hypocoercivity}


In the present section, we wish to study the decay generated by the linearized operator (\ref{E:linear}).  In fact, we will
study an approximation to that operator obtained by dropping the non-local operator $(\sin x) \partial_y \Delta^{-1}$, and
studying the approximate equation
\begin{eqnarray} \nonumber
\partial_t w &=& \cLa w \\ \label{E:linear_approx}
\cLa w &=& \nu \Delta w - a e^{-\nu t} (\sin x) \partial_y w\ .
\end{eqnarray}
Note that the term we dropped is a relatively compact perturbation of $\cLa$ and hence we believe that $\cLa$
should capture the most important aspects of the linearized evolution.  Furthermore
we present numerical results on the full linearized operator in Section \ref{S:numerics}
that support this point of view.  We discuss at the end of this section the extension of our current
results to the full linear operator.  We also note that a similar approach has recently been utilized
to study the linearization of the two-dimensional Navier-Stokes equation around the Oseen vortex \cite{Deng12}.

To state our results, note that if we write the solution $w  = \sum_{k,\ell } \hat{w}(k,\ell) e^{i(kx+\ell y)}$ of \eqref{E:linear_approx} as a Fourier series, the 
operator $\cLa$ does not mix different values of the Fourier index $\ell$.  Thus, if we define
\begin{equation}
w(x,y) = \sum_{\ell} (\hat{w}_{\ell}(x) ) e^{i \ell y} = \sum_{\ell} \left( \sum_k \hat{w}(k,\ell) e^{ikx} \right) e^{i \ell y}
\end{equation}
we can rewrite
\begin{equation}
\cLa w = \sum_{\ell} (\cLaell \hat{w}_{\ell} ) e^{i \ell y}
\end{equation}
where
\begin{equation}\label{E:linear_ell}
\cLaell \hat{w}_{\ell}  = \nu (\partial_x^2 - \ell^2) \hat{w}_{\ell} - i a \ell e^{-\nu t} (\sin x) \hwell\ .
\end{equation}

Motivated by Villani's work on hypocoercivity, \cite{Villani09}, we define operators,
\begin{eqnarray} \label{E:BCdef}
\Bell \hwell & = & -i a \ell e^{-\nu t} (\sin x) \hwell \\ \nonumber
\Cell \hwell &=& [\partial_x, \Bell ] \hwell = -i a \ell e^{-\nu t} (\cos x) \hwell\ .
\end{eqnarray}
As we will see, the fact that $ [\partial_x, \Bell ] \ne 0$ is responsible for the faster than expected decay of
solutions of our linear equation.

\begin{Remark}  The discussion of anomalous modes from the previous section is simplified for $\cLa$.
The operators $\Bell$ and $\Cell$ both annihilate any functions $\hat{w}_0$, while there are no nonzero elements of
the kernel of $\Bell$ or $\Cell$ for any $\ell \ne 0$.  Thus, we will show that the only anomalous modes are
those with $\ell =0 $ and we have accelerated decay for all solutions of \eqref{E:linear_approx} of the form
\begin{equation}
w(x,y) = \sum_{\ell \ne 0} \hwell (x) e^{i \ell y}\ .
\end{equation}
\end{Remark}

Given the preceding remark, we define the Banach space
\begin{equation} \label{E:def_X}
X = \left\{ w \in L^2(\mathbb{T}^2): \hat{w}_0 = 0 ,\ \ \sum_{l \neq 0} \left[\| \hwell \|^2 + \sqrt{\frac{\nu}{| \ell |}} \| \partial_x \hwell  \|^2 +  \frac{1}{\sqrt{\nu}| \ell |^{3/2}} \|C^\ell\hwell \|^2\right]. =: \|w\|_X^2 < \infty \right\}\ .
\end{equation}
 Note that the $l = 0$ modes are absent from functions in $X$.
 In this definition, and throughout the remainder of the paper, $\| \cdot \| $ will denote the norm on $L^2(\mathbb{T}^1)$.
We will use $(\cdot, \cdot)$ to denote the innerproduct on $L^2(\mathbb{T}^1)$.

\begin{Theorem} \label{Thm:bar_conv} Given any constant $\tau$ and $T  \in [0, \tau/\nu]$, there exist constants $K$ and $M$ that are $\mathcal{O}(1)$ with respect to $\nu$ such that the following holds. If $\nu$ is sufficiently small, then the solution to the linear equation (\ref{E:linear_approx}) with initial condition $w^0 \in X$ satisfies the estimate
\[
\|w(t)\|_{X}^2 \leq K e^{-M\sqrt{\nu}t} \|w^0\|_X^2
\]
for all $t \in [0, T]$. 
\end{Theorem}

Hence, at least in our approximation, the bar states are more stable than one would expect from the diffusive decay of $\mathcal{O}(e^{- \nu t})$. In fact, if we take $T = \tau/\nu$, we find $e^{-M\sqrt{\nu}T} = e^{-M\tau } \ll e^{-\nu T} = e^{-\tau }$. (See Remark \ref{R:M_0}, below, for a discussion of how the constants $M$ and $K$ are related to $T$.) Therefore, even though this result is only valid for finite times, it still implies that, for initial data in the above mentioned subspace, there will be an initial period of rapid decay to the bar states, which then themselves decay on the diffusive timescale.

As we mentioned in the introduction, this result is reminiscent of our prior work on metastable states in Burgers equation,
\cite{BeckWayne11}.  As in that work, the bar states form an invariant manifold (in fact, an invariant subspace in this
case) which connects to the long-time asymptotic state of the system (the zero state in this case), and at that point, it is tangent
to the eigenspace of one of the eigenvalues with smallest non-zero value.  Furthermore, the preceding theorem shows that
(at least in our approximation), nearby solutions are attracted to the invariant manifold at a rate much faster than the motion
along the manifold, justifying the identification of these solutions as metastable states of the system.

The remainder of this section is devoted to the proof of this Theorem. 

Consider equation (\ref{E:linear_approx}).   As we said above, we wish to use Villani's hypercoercivity method,
as developed by  \cite{GallagherGallayNier09}, and as a consequence, the properties of the operators $\Bell$ and
$\Cell$ introduced in the introduction of this section will play an important role.  Thus, we begin, with a simple
lemma detailing some of the properties of these operators.

\begin{Lemma} The operators $\Bell$ and $\Cell$ are bounded, anti-symmetric operators on $L^2(\mathbb{T}^1)$
for any $\ell \ne 0$, with bound
\begin{equation}
\| \Bell \hwell \|  \le | a  \ell | e^{-\nu t} \| \hwell \| ;\ \ \| \Cell \hwell \|  \le | a  \ell | e^{-\nu t} \| \hwell \|\ .
\end{equation}
Furthermore,  $\Bell \Cell = \Cell \Bell$ - i.e. these operators commute.
\end{Lemma}

\begin{Proof}  The proof is a simple computation using the definition of $\Bell$ and $\Cell$.  For instance
\begin{equation*}
\| \Bell \hwell \|^2 = (a \ell)^2 e^{-2\nu t} \int \sin^2 x | \hwell (x) |^2 dx \le (a \ell)^2 e^{-2\nu t} \| \hwell \|^2
\end{equation*}
while 
\begin{equation*}
(\hwell, (\Bell  \hat{u}_{\ell}) ) = a \ell e^{-\nu t} \int \overline{\hwell (x)} (-i  \hat{u}_{\ell}(x) ) dx = - a \ell e^{-\nu t}
 \int (\overline{- i \hwell (x) }) \hat{u}_{\ell}(x) dx = - ( \Bell \hwell,  \hat{u}_{\ell})
\end{equation*}

with a similar computation for $\| \Cell \hwell \|$.\end{Proof}

Because we wish to follow the scheme developed by Gallagher, Gallay and Nier, we 
we rescale the time by $t \to \nu t$ and consider the following problem equivalent to \eqref{E:linear_ell}, 
\begin{equation} \label{eq:linear_eq_two}
\partial_t \hwell  =    (\partial_x^2-\ell^2) \hwell   + \frac{1}{\nu}  \Bell \hwell,\ \ \ell \ne 0\ .
\end{equation}
Following \cite{Villani09} and  \cite{GallagherGallayNier09}, we will measure the changes in solutions of \eqref{eq:linear_eq_two}
in terms of the functionals 
\begin{equation}
\Phi^{\ell} (t) = \| \hwell \|^2 +  \alpha \| \partial_x \hwell \|^2 - 2\beta  \Re (\partial_x \hwell , \Cell  \hwell)+ \gamma (\Cell \hwell,\Cell \hwell).
\end{equation}
We will establish decay for solutions of \eqref{eq:linear_eq_two} for each $\ell$ by choosing the constants
$\alpha$, $\beta$ and $\gamma$ appropriately, and then combine these estimates to prove Theorem \ref{Thm:bar_conv}.

\begin{Remark}\label{rem:norm} Note that 
\[
| 2 \beta \Re (\partial_x \hwell, \Cell  \hwell) | \le \frac{\alpha}{2} \|( \hwell)_x \|^2 + \frac{2\beta^2}{\alpha} \| \Cell \hwell\|^2 \ .
\]
Thus, if 
\begin{equation}\label{eq:abc_relationship}
\beta^2 < \frac{ \alpha \gamma}{4},
\end{equation}
we have
\[
\| \hwell \|^2 + \frac{\alpha}{2} \| (\hwell)_x \|^2 + \frac{\gamma}{2} \| \Cell  \hwell \|^2 < \Phi^{\ell} (t) < 
\| \hwell \|^2 + \frac{3 \alpha}{2} \| (\hwell)_x \|^2 + \frac{3 \gamma}{2} \| \Cell \hwell \|^2.
\]
\end{Remark}

Consider the time change of this functional.  Also, for the moment, we will be working with equation \eqref{eq:linear_eq_two}
for a fixed value of $\ell$, so we suppress the subscript $\ell$ on the functions $\hwell$ to avoid overburdening the notation.
\begin{eqnarray}\label{eq:dt} \nonumber
&& \frac{d}{dt} \Phi^{\ell} (t) =  \left( ( \hw_t, \hw)+(\hw,\hw_t) \right) +
\alpha \left( ( \partial_x \hw_t, \partial_x \hw) + (\partial_x \hw, \partial_x \hw_t) \right) \\ 
&& \qquad \qquad 
 -2\beta \Re \left( (\partial_x \hw_t , \Cell \hw)+ (\partial_x \hw , \Cell \hw_t) \right)
 + \gamma \left( (\Cell \hw_t,\Cell\hw) + (\Cell\hw,\Cell \hw_t) \right) \\ \nonumber
 && \qquad \qquad   -2\beta \Re (\partial_x \hw , \frac{ d\Cell}{dt} \hw ) + \gamma \left(   
 ( \frac{ d\Cell}{dt}   \hw,\Cell\hw) + (\Cell\hw,\frac{ d\Cell}{dt} \hw)  \right) .
\end{eqnarray}
We bound in turn each of the pairs of terms in \eqref{eq:dt}


\subsection{The term $ \left( (\hw_t, \hw)+(\hw,\hw_t) \right)$}

Using \eqref{eq:linear_eq_two} we can rewrite this as:
\begin{equation}
\left( (  (-\ell^2 + \partial_x^2 + \frac{1}{\nu} \Bell)\hw, \hw ) + (\hw,  (-\ell^2 + \partial_x^2 + \frac{1}{\nu  } \Bell) \hw) \right)
= -2 \ell^2 \| \hw \|^2  -2  \| \hw_x \|^2, \label{E:term1}
\end{equation}
where we used the anti-symmetry of $\Bell$ to eliminate terms like $(\hw,\Bell\hw)$.


\subsection{The $\alpha$ term $\left( (\hw_{xt}, \hw_x)+(\hw_x,\hw_{xt}) \right) $}

Again, we rewrite this as:
\begin{eqnarray*} \nonumber
\left( ( \partial_x  (-\ell^2 + \partial_x^2 + \frac{1}{\nu} \Bell)\hw, \hw_x) + 
(\hw_x,  \partial_x (-\ell^2 + \partial_x^2 + \frac{1}{\nu} \Bell)\hw \right)
&=& -2 \ell^2 \| \hw_x \|^2 - 2 \| \hw_{xx} \|^2\  \\ 
&&\quad  +\frac{1}{\nu} [ (\partial_x (\Bell\hw), \hw_x)+(\hw_x,\partial_x (\Bell\hw))]
\end{eqnarray*}
Using the anti-symmetry of $\Bell$, we can rewrite the last two terms as
\begin{eqnarray*}
[ (\partial_x (\Bell\hw), \hw_x)+(\hw_x,\partial_x (\Bell\hw))]  &=& 
(\Bell \hw_x,\hw_x) + ([\partial_x,\Bell] \hw,\hw_x) +(\hw_x,\Bell \hw_x) +(\hw_x,[\partial_x,\Bell] \hw)\\
&=& 2\Re (\hw_x,\Cell\hw)\ .
\end{eqnarray*}

Since $|2\Re (\hw_x,\Cell\hw)| \le 2\| \hw_x \| \|\Cell \hw\|$, we finally obtain
\begin{equation} 
\alpha \left( (-\ell^2 + \partial_x^2 + \frac{1}{\nu} \Bell)\hw, \hw_x) + 
(\hw_x,  \partial_x (-\ell^2 + \partial_x^2 + \frac{1}{\nu} \Bell) \hw\right)
\le  -2 \alpha\ell^2 \| \hw_x \|^2 - 2\alpha \| \hw_{xx} \|^2 +    \frac{2\alpha}{\nu}\| \hw_x \|\  \|\Cell \hw\|.
\label{E:term2}
\end{equation}


\subsection{The $\beta$ term $\left( (\partial_x \hw_t , \Cell \hw)+ (\partial_x \hw , \Cell \hw_t) \right)$}

Using \eqref{eq:linear_eq_two} we can rewrite this as:
\begin{eqnarray*}
(\partial_x \hw_t , \Cell \hw)+ (\partial_x \hw , \Cell \hw_t) &=& (\partial_x(-\ell^2 + \partial_x^2 + \frac{1}{\nu} \Bell )\hw, \Cell\hw) + (\partial_x \hw, \Cell(-\ell^2 + \partial_x^2 + \frac{1}{\nu} \Bell ) \hw ) \\ 
&=& -2 \ell^2 \mathrm{Re}(\partial_x \hw, \Cell \hw) + [(\hw_{xxx},\Cell\hw) + (\hw_x,\Cell \hw_{xx})] \\
&&\qquad + \frac{1}{\nu}[ (\partial_x (\Bell \hw),\Cell\hw) + (\hw_x, \Cell(\Bell \hw)) ]
\end{eqnarray*}

By Cauchy-Schwartz inequality, we can bound
$
| (\hw_x,\Cell \hw)| \le  \| \hw_x \|\  \| \Cell \hw\|\ .
$
Integrating by parts, and using the fact that $[\partial_x, \Cell] = -\Bell $, we can write
\begin{eqnarray*}
(\hw_{xxx},\Cell \hw) + (\hw_x,\Cell \hw_{xx}) &=& -(\hw_{xx} , \partial_x (\Cell \hw)) + (\hw_x,\Cell \hw_{xx}) \\ 
&=& -2 \Re(\hw_{xx} ,\Cell \hw_x ) + (\hw_{xx},\Bell  \hw) \\
&=& -2 \Re(\hw_{xx} ,\Cell \hw_x ) - (\hw_x, \partial_x \Bell  \hw) \\
&=& -2 \Re(\hw_{xx} ,\Cell \hw_x ) - (\hw_x,  \Bell  \hw_x)- (\hw_x, \Cell  \hw) \ .
\end{eqnarray*}
Because $\Bell $ is anti-symmetric, the middle term has zero real part, and so it will not contribute anything
 to the time derivative of $\Phi^{\ell}(t)$, as we only consider the real part of this term.  We can therefore bound this expression by
\begin{equation*}
| (\hw_{xxx},\Cell \hw) + (\hw_x,\Cell  \hw_{xx})| \le  2  \| \hw_{xx} \| \| \Cell  \hw_x \| + \| \hw_x \| \|\Cell  \hw\|
\end{equation*}
Next, commuting $\partial_x$ and $\Bell $, we have
\begin{eqnarray*}
(\partial_x (\Bell \hw),\Cell \hw) + (\hw_x, \Cell (\Bell \hw)) &=& (\Bell  \hw_x, \Cell \hw) 
+ (\Cell \hw,\Cell \hw) + (\hw_x, \Bell  \Cell \hw) + ( \hw_x, [\Cell ,\Bell ] \hw)  \\ 
&=& \| \Cell  \hw \|^2 + (\hw_x, [\Cell ,\Bell ]\hw) = \| \Cell  \hw \|^2\ ,
\end{eqnarray*}
where the penultimate equality used the anti-symmetry of $\Bell $ and the last equality used the fact that $[\Bell, \Cell]=0$.

Thus, combining the above estimates, we obtain
\begin{equation}\label{E:term3}
-2 \beta \mathrm{Re} [(\partial_x \hw _t , \Cell  \hw )+ (\partial_x \hw  , \Cell  \hw _t)]  \leq 
(4\ell^2 +2)\beta \|\hw _x\| \|\Cell \hw \| + 4\beta\|\hw _{xx}\|\|\Cell \hw _x\| - \frac{2\beta}{\nu}\|\Cell \hw \|^2.
\end{equation}


\subsection{The $\gamma$ term $ \left( (\Cell  \hw _t,\Cell \hw ) + (\Cell \hw ,\Cell  \hw _t) \right)$}

Rewrite this as
\begin{eqnarray*}
&&( \Cell  (-\ell^2 + \partial_x^2 + \frac{1}{\nu} \Bell ) \hw , \Cell \hw ) + 
(\Cell \hw , \Cell  (-\ell^2 + \partial_x^2 + \frac{1}{\nu} \Bell ) \hw ) \\
&&\qquad \qquad = - 2 \ell^2 \|\Cell \hw \|^2 + \left( (\Cell  \hw _{xx}, \Cell  \hw )
 + (\Cell \hw  , \Cell  \hw _{xx} ) \right) + \frac{1}{\nu} \left( (\Cell  (\Bell \hw ), \Cell \hw ) + (\Cell \hw , \Cell  (\Bell \hw ))  \right)
\end{eqnarray*}
We have
\begin{eqnarray*}
(\Cell  \hw _{xx}, \Cell  \hw ) + (\Cell \hw  , \Cell  \hw _{xx} ) &=& - [ ( \Cell  \hw _x, \partial_x (\Cell  \hw ) ) + (\partial_x( \Cell  \hw ), \Cell  \hw _x)] + 
[(\Bell \hw _x, \Cell \hw )+(\Cell  \hw , \Bell  \hw _x)] \\ 
&=& - [ ( \Cell  \hw _x, \Cell \partial_x \hw  - \Bell  \hw ) + (\Cell \partial_x \hw  - \Bell \hw , \Cell  \hw _x)] \\
&&\qquad + [(\Bell \hw _x, \Cell \hw )+(\Cell  \hw , \Bell  \hw _x)] \\ 
&=& -2 ( \Cell  \hw _x, \Cell  \hw _x) + (\partial_x \Cell \hw  + \Bell \hw , \Bell \hw ) + (\Bell \hw ,  \partial_x \Cell \hw  + \Bell \hw ) \\
&& \qquad + [(\Bell \hw _x, \Cell \hw )+(\Cell  \hw , \Bell  \hw _x)] \\
&=& -2 ( \Cell  \hw _x, \Cell  \hw _x) + 2(\Bell \hw , \Bell \hw ) - [(\Cell \hw , \partial_x \Bell \hw ) + (\partial_x \Bell \hw , \Cell \hw )] \\
&&\qquad + [(\Bell \hw _x, \Cell \hw )+(\Cell  \hw , \Bell  \hw _x)] \\
&=& -2 ( \Cell  \hw _x, \Cell  \hw _x) - 2 ( \Cell  \hw , \Cell  \hw ) + 2(\Bell \hw , \Bell \hw ).
\end{eqnarray*}
The  term proportional to $\nu^{-1}$ can be computed using the fact that
\[
(\Cell  (\Bell \hw ), \Cell \hw ) + (\Cell \hw , \Cell  (\Bell \hw )) = (\Bell  (\Cell \hw ), \Cell \hw ) + (\Cell \hw , \Bell  (\Cell \hw ))= 0,
\]
where the first step used the fact that $\Cell$ and $\Bell$ commute, and
the second step used the anti-symmetry of $\Bell$.

Hence, we have
\begin{equation} \label{E:term4}
\gamma[(\Cell  \hw_t,\Cell \hw) + (\Cell \hw,\Cell  \hw_t)] \leq 
 -(2\ell^2+2)\gamma \|\Cell \hw\|^2 -2\gamma \|\Cell  \hw_x\|^2 + 2\gamma\|\Bell \hw\|^2.
\end{equation}

\subsection{The terms containing $\frac{d \Cell}{dt}$}

Noting that $\frac{d \Cell}{dt} = -\nu \Cell$, (from \eqref{E:BCdef}), we have
\begin{equation}
|2\beta \Re (\partial_x \hw , \frac{ d\Cell}{dt} \hw )| \le 2 \beta \nu \| \hw_x\| \| \Cell \hw \|
\end{equation}
while
\begin{equation}
| \gamma \left( ( \frac{ d\Cell}{dt}   \hw , \Cell \hw ) + (\Cell \hw ,  \frac{ d\Cell}{dt}   \hw ) \right) | = 2 \nu \gamma \| \Cell \hw \|^2
\end{equation}

\subsection{Estimate on the time evolution of $\Phi^{\ell}$}
Combining equations (\ref{E:term1})-(\ref{E:term4}), we obtain
\begin{eqnarray*}
\frac{d}{dt}\Phi^{\ell}(t) &\leq& 
-2 \ell^2\| \hw \|^2 - [ 2+ 2\alpha \ell^2] \|\hw_x\|^2 - 2\alpha \|\hw_{xx}\|^2 + 
[ \frac{2\alpha}{\nu} + 2\beta (2\ell^2 + 1+\nu )] \|\hw_x\| \|\Cell\hw\| \\
&&\quad + 4\beta \|\hw_{xx}\|\|\Cell\hw_x\| -[ (2\ell^2+2) \gamma + \frac{2\beta}{\nu}-2 \gamma \nu ] 
\|\Cell\hw\|^2 -2\gamma  \|\Cell \hw_x\|^2  + 2\gamma \|\Bell\hw\|^2.
\end{eqnarray*}
We now use the facts that
\begin{eqnarray*}
4\beta \|\hw_{xx}\|\|\Cell\hw_x\| &\leq& \alpha \|\hw_{xx}\|^2 + \frac{4 \beta^2}{\alpha} \|\Cell\hw_x\|^2 \\
\frac{2\alpha}{\nu}\|\hw_x\|\|\Cell\hw\| &\leq&  \|\hw_x\|^2 + \frac{\alpha^2}{\nu^2}\|\Cell\hw\|^2 \\ 
2\beta(2\ell^2 + 1)\|\hw_x\|\|\Cell\hw\| &\leq& \frac{1}{2}\alpha(2\ell^2 + 1+\nu)\|\hw_x^2\| + \frac{2\beta^2(2\ell^2+1+\nu)}{\alpha} 
\|\Cell\hw\|^2
\end{eqnarray*}
to obtain
\begin{eqnarray*}
\frac{d}{dt}\Phi^{\ell}(t) &\leq& -2 \ell^2 \| \hw \|^2 + 2\gamma \|\Bell\hw\|^2 - [ 2+ 2\alpha \ell^2 - 1 - \frac{1}{2}\alpha(2\ell^2 + 1+\nu)] \|\hw_x\|^2 - (2\alpha -\alpha) \|\hw_{xx}\|^2 \\
&&\qquad - [2\gamma - \frac{4\beta^2}{\alpha}] \|\Cell\hw_x\|^2 - [\frac{2\beta}{\nu} + (2\ell^2 + 2)\gamma - 
\frac{\alpha^2}{\nu^2} - \frac{2 \beta^2(2\ell^2 + 1+\nu )}{\alpha}-2 \gamma \nu] \|\Cell\hw\|^2.
\end{eqnarray*}
Now set $\alpha = \alpha_0 \sqrt{\nu}$, $\beta = \beta_0$, and $\gamma = \gamma_0/\sqrt{\nu}$, where $\alpha_0$, 
$\beta_0$, and $\gamma_0$ will be assumed to be independent of $\nu$,
to find
\begin{eqnarray*}
\frac{d}{dt}\Phi^{\ell}(t) &\leq& -2 \ell^2 \| \hw \|^2 + 2\frac{\gamma_0}{\sqrt{\nu}} \|\Bell\hw\|^2 - [ (1  -\nu  - \frac{\alpha_0\sqrt{\nu}}{2}) + \alpha_0\sqrt{\nu}\ell^2] \|\hw_x\|^2 - \alpha_0 \sqrt{\nu} \|\hw_{xx}\|^2 \\
&&\qquad - \frac{2(\gamma_0\alpha_0 - 2\beta_0^2)}{\alpha_0 \sqrt{\nu}} \|\Cell\hw_x\|^2 - [\frac{(2\beta_0 - \alpha_0^2)}{\nu} + 2\ell^2 \frac{(\gamma_0\alpha_0 - 2\beta_0^2)}{\alpha_0\sqrt{\nu}}  \\ \nonumber && \qquad \qquad + 
2\frac{(\gamma_0 \alpha_0 - \beta_0^2)}{\sqrt{\nu}\alpha_0}-\frac{2 \beta_0^2}{\alpha_0} \sqrt{\nu} - 2 \gamma_0 \sqrt{\nu} ] 
\|\Cell\hw\|^2.
\end{eqnarray*}
Hence, as long as $\alpha_0$, $\beta_0$, and $\gamma_0$ are chosen so that (\ref{eq:abc_relationship}), and
\begin{equation} \label{eq:ab_condition}
\beta_0 \geq 4 \alpha_0^2 
\end{equation}
hold,  and if $\nu$ is sufficiently small, we have
\begin{equation}\label{eq:upper_bound_one}
\frac{d}{dt}\Phi^{\ell}(t) \le -2 \ell^2 \| \hw \|^2 + 2\frac{\gamma_0}{\sqrt{\nu}}
 \|\Bell\hw\|^2 - \frac{1}{4} \|\hw_x\|^2 - \alpha_0 \sqrt{\nu} \|\hw_{xx}\|^2 - 
 \frac{\gamma_0}{\sqrt{\nu}} \|\Cell\hw_x\|^2 - \frac{3\beta_0}{2\nu}\|\Cell\hw\|^2.
\end{equation}

Our goal is to show that $(\Phi^{\ell})^{\prime}(t) \le \frac{M}{\sqrt{\nu}} \Phi^{\ell}(t)$.  All terms in \eqref{eq:upper_bound_one} are of the correct form except for $-2\ell^2\| \hw \|^2 + \frac{2 \gamma_0}{\sqrt{\nu}}\|\Bell\hw\|^2$.  We handle these terms with the aid of a proposition modeled after the analogous bound of \cite[Proof of Prop 4.1]{GallagherGallayNier09}.

\begin{Proposition} \label{prop:harmonic}
If $\tau$, $T \in [0, \tau/\nu]$, $a$, and $| \ell | > 1$ are fixed, there exists a constant $M_0$ such that, 
for $\nu$ sufficiently small and for all $t \in [0, T]$, 
\begin{equation}\label{E:mod_prop}
\frac{1}{8} \| u_x \|^2 + \frac{ \beta_0 }{2 \nu} \| \Cell u \|^2 \ge \frac{M_0 | \ell |\sqrt{\beta_0} }{\sqrt{\nu}}\| u \|^2,
\end{equation}
for any $u \in H^1(\mathbb{T}^1)$.
\end{Proposition}

\begin{Remark}\label{R:M_0}
The constant $M_0$ can be chosen uniformly
for $t \in [0,T]$, but its value will depend on $T$.
The  proof  below indicates that $M_0 \sim  \mathcal{O}(e^{-\nu T})$. 
This makes sense: since $\| \Cell \| \to 0$ as $t \to \infty$, we must have 
$M_0 \to 0$ as $t \to \infty$, as well.  This implies we may take $T = \mathcal{O}(1/\nu)$, 
since then $M_0 = \mathcal{O}(1)$, in the statement of this proposition, and therefore also in the statement of Theorem \ref{Thm:bar_conv}. 
\end{Remark}

Before proving this proposition, let's look at its consequences  for equation \eqref{eq:upper_bound_one}. 
Since $\| \Bell \hw \| \le a le^{-\nu t} \| \hw \| \le al \| \hw \|$, equation (\ref{E:mod_prop}) implies 
\begin{equation}\label{eq:l2bound}
\frac{2 \gamma_0}{\sqrt{\nu}} \| \Bell\hw \|^2
 - \left( \frac{1}{8} \|  \hw_x \|^2 + \frac{\beta_0}{2\nu} \| \Cell \hw \|^2 \right) 
\le - (\frac{M_0  l\sqrt{\beta_0}}{\sqrt{\nu}}- \frac{2\gamma_0 a^2l^2}{\sqrt{\nu}} ) \| \hw \|^2.
\end{equation}

We now define 
\begin{equation}\label{eq:b_def}
\beta_0 = \frac{16 \gamma_0^2 a^4\ell^2}{M_0^2},
\end{equation}
which implies that the right hand side of \eqref{eq:l2bound} is less than or equal to
\[
-\frac{2 \gamma_0a^2\ell^2}{\sqrt{\nu}} \| \hw \|^2.
\]
Next, to insure that \eqref{eq:ab_condition} holds, we set
\begin{equation}\label{eq:a_def}
\alpha_0^2 = \frac{\beta_0}{4} = \frac{4 \gamma_0^2 a^4 \ell^2}{M_0^2}.
\end{equation}
Finally, we must check that \eqref{eq:abc_relationship} is satisfied.  Given the definitions of $\alpha_0$ and $\beta_0$, this inequality will hold provided
\[
\gamma_0 > \frac{4 \beta_0^2}{\alpha_0} = 8 \beta_0^{3/2} = 8 \left(\frac{16 \gamma_0^2 a^4 \ell^2}{M_0^2 }\right)^{3/2}.
\]
Therefore, combining this with equation (\ref{eq:b_def}) and (\ref{eq:a_def}), we take
\begin{equation}\label{E:def_g}
\gamma_0 = \frac{M_0^{3/2}}{64\sqrt{2}a^3 |\ell|^{3/2}}, \qquad \alpha_0 = \frac{M_0^{1/2}}{32\sqrt{2} a |\ell|^{1/2}}, \qquad \beta_0 = \frac{M_0}{32(16) a^2 |\ell|}.
\end{equation}
With these definitions, we obtain:
\begin{equation} \label{E:phidot_est}
\frac{d}{dt} \Phi^{\ell} (t)  \le    -2\ell^2 \|\hw\|^2 - \frac{2 \gamma_0 a^2 l^2 }{\sqrt{\nu}}  \| \hw \|^2 - \frac{1}{8} \| \hw_x \|^2 - \alpha_0 \sqrt{ \nu}  \| \hw_{xx} \|^2 
- \frac{\beta_0 }{ \nu} \| \Cell  \hw \|^2 - \frac{\gamma_0 }{\sqrt{\nu}} \| \Cell  \hw_x \|^2 
\end{equation}
where $\alpha_0$, $\beta_0$, and $\gamma_0$ are defined in \eqref{E:def_g}.  

We wish to show that $(\Phi^{\ell})^{\prime}  \leq -(M/\sqrt{\nu})\Phi$ for some $M$, which will hold if $\Phi/\sqrt{\nu}$ is less than or equal to $-1/M$ times the right hand side of \eqref{E:phidot_est}. Remark \ref{rem:norm} implies that
\[
\frac{1}{\sqrt{\nu}} \Phi^{\ell}  \leq \frac{1}{\sqrt{\nu}}\|\hw\|^2 + \frac{3\alpha_0}{2} \|\hw_x\|^2 + \frac{3\gamma_0}{2\nu} \|\Cell \hw\|^2,
\] 
which needs to be less than or equal to 
\[
\frac{1}{M} \left( 2 \ell^2 \|\hw\| + \frac{2 \gamma_0 a^2 l^2 }{\sqrt{\nu}}  \| \hw \|^2 
+ \frac{1}{8} \| \hw_x \|^2 + \alpha_0 \sqrt{ \nu}  \| \hw_{xx} \|^2 
+ \frac{\beta_0 }{ \nu} \| \Cell  \hw \|^2 + \frac{\gamma_0 }{\sqrt{\nu}} \| \Cell  \hw_x \|^2 
 \right).
\]
Using equation \eqref{E:def_g}, we find this will hold if we chose $M$ such that
\begin{eqnarray*}
M &\leq & 2\gamma_0 a^2 \ell^2 = \frac{M_0^{3/2}\sqrt{| \ell |}}{32\sqrt{2}a},\\
M &\leq& \frac{1}{12 \alpha_0} = \frac{8\sqrt{2} a \sqrt{| \ell |}}{3\sqrt{M_0}} \\
M &\leq& \frac{2\beta_0}{3\gamma_0} = \frac{24 a \sqrt{| \ell |}}{\sqrt{2M_0}}. 
\end{eqnarray*}
Since $| \ell | \geq 1$ on the subspace $\mathcal{M}$, we can take $M= \mathcal{O}(M_0^{3/2})$, independent of $\ell$. Although the proof of Proposition \ref{prop:harmonic}, below, implies that $M = \mathcal{O}(M_0^{3/2}) = \mathcal{O}(e^{-\frac{3\nu}{2}t})$, this will still be an $\mathcal{O}(1)$ quantity with respect to $\nu$ for all $t \leq \mathcal{O}(1/\nu)$. Thus, we may take $T = \mathcal{O}(1/\nu)$ in the statement of the Proposition. Gronwall's inequality then implies that
\begin{equation}\label{E:phi_decay}
\Phi^{\ell} (t) \leq \Phi^{\ell} (0) e^{-\frac{M}{\sqrt{\nu}} t}.
\end{equation}

We now use this estimate to control the norm of solutions. Throughout this estimate, $C_1, C_2, K$ denote generic constants that are independent of any relevant quantity. Equation \eqref{E:def_g} and Remark \ref{rem:norm} imply that
\begin{equation}\label{E:pinch}
\| \hw \|^2 + \frac{1}{2} C_1 \sqrt{\frac{\nu}{|\ell |}} \| \hw_x \|^2 + \frac{1}{2}  C_2 \frac{1}{\sqrt{\nu}|\ell |^{3/2}} \| C^\ell \hw \|^2 < \Phi^{\ell} (t) < 
 \| \hw \|^2 + \frac{3}{2}  C_1 \sqrt{\frac{\nu}{| \ell |}} \| \hw_x \|^2 + \frac{3}{2} C_2 \frac{1}{\sqrt{\nu}| \ell |^{3/2}} \| C^\ell \hw \|^2.
\end{equation}
Recall that
\[
\|w\|_X^2 = \sum_{l \neq 0} \left[\| \hwell \|^2 + \sqrt{\frac{\nu}{| \ell |}} \| \partial_x \hwell  \|^2 +  \frac{1}{\sqrt{\nu}| \ell |^{3/2}} \|C\hwell \|^2\right].
\]
Then, \eqref{E:pinch} and \eqref{E:phi_decay} imply
\[
\|w(t)\|_X^2 \leq K e^{-M\sqrt{\nu}t} \|w(0)\|_X^2,
\]
and we obtain Theorem \ref{Thm:bar_conv}.


\begin{Proof}[ of Proposition \ref{prop:harmonic}]
We now prove the proposition via a series of lemmas, following closely the approach of Gallagher, Gallay and Nier.


Define a partition of unity $\chi^2_{+}(x) + \chi^2_{-}(x) = 1$ where $\chi_{\pm}$ are $C^{\infty}$
and non-negative with $supp(\chi_{+}) \subset [-\delta, \pi+\delta]$, for some $\delta$ small and positive, and$\chi_{+}(x) =1 $ for all $x\in [\delta,\pi-\delta]$, $supp(\chi_{-}) \subset [-\pi-\delta, \delta]$,
and $\chi_{-}(x) =1 $ for all $x\in [-\pi+\delta, -\delta]$. For any $u \in L^2$, define $u_{\pm} = \chi_{\pm} u$.

\begin{Lemma} \label{lem:plus_minus} There exists a constant $C>0$ such that for any $u\in H^1$, one has
\[
\| \partial_x u \|^2 \ge \frac{2}{3} (\| \partial_x u_{+} \|^2 + \| \partial_x u_{-} \|^2) - C \| u \|^2
\]
\[
(\| \partial_x u_{+} \|^2 + \| \partial_x u_{-} \|^2) \ge \frac{3}{4} \| \partial_x u\|^2 - C \| u \|^2\ .
\]
\end{Lemma}
\proof First note that
\begin{eqnarray}\label{eq:IML}
\int | \partial_x u |^2 &=& \int \partial_x \overline{u} (\chi^2_{+} + \chi^2_{-}) \partial_x u
= \int \partial_x \overline{u} \chi^2_{+} \partial_x u + \int \partial_x \overline{u} \chi^2_{-} \partial_x u. 
\end{eqnarray}
One can write 
\begin{eqnarray*}
\int \partial_x \overline{u} \chi^2_{+} \partial_x u &=& \int (\chi_{+} \partial_x \overline{u})  (\chi_{+} \partial_x u)
\\ \nonumber
& =&  \int ( \partial_x \overline{u_{+}})( \partial_x u_{+}) -\left( \int ((\partial_x \chi_{+} )\overline{u})(\chi_{+} \partial_x  u) + \int ( \partial_x \overline{u_{+}}) ( (\partial_x \chi_{+}) u )\right).
\end{eqnarray*}
We can now bound
\[
| \int ( \partial_x \overline{u_{+}}) ( (\partial_x \chi_{+}) u )| \le C_1  \| u \| \| \partial_x u \|
\le \frac{1}{4} \| \partial_x u \|^2 + C_2 \| u \|^2.
\]
Applying the same sort of estimates to the last term on the right hand side of \eqref{eq:IML}, we obtain
\[
\int | \partial_x u |^2 \ge \left( \| \partial_x u_{+} \|^2 + \| \partial_x u_{-}\|^2 \right) - \| \partial_x u \|^2
- 4C_2 \| u \|^2,
\]
from which the first inequality in the lemma follows. 

To prove the second inequality, we again use equation (\ref{eq:IML}). Consider the first of the two integrals on the right hand side, which can be written
\begin{eqnarray*}
\int \overline{(\chi_{+} u_x)} (\chi_{+} u_x) &=& \int \overline{( \partial_x u_{+} - \chi_{+}' u)}( \partial_x u_{+} - \chi_{+}' u) \\ \nonumber
&=& \int | \partial_x u_{+} |^2 + \int (\chi_{+}')^2 |u|^2 - 2 \Re \int (\overline{\partial_x u_{+}} )\chi_{+}' u \\ \nonumber
&\le & \frac{4}{3}  \int | \partial_x u_{+} |^2  + C_3  \int  |u|^2 ,
\end{eqnarray*}
where the last inequality used the bounds
 $\int (\chi_{+}')^2 |u|^2 \le C_1 \| u \|^2$ and $|2 \Re \int (\overline{\partial_x u_{+}} )\chi_{+}' u|
\le \frac{1}{3} \int | \partial_x u_{+} |^2  + C_2  \int  |u|^2$. If we apply a similar estimate to the term $\int \overline{(\chi_{-} u_x)} (\chi_{-} u_x)$ in (\ref{eq:IML}),
we then obtain 
\[
\|u_x\|^2 = \int \partial_x \overline{u} \chi^2_{+} \partial_x u + \int \partial_x \overline{u} \chi^2_{-} \partial_x u \leq \frac{4}{3} \left(   \int | \partial_x u_{+} |^2 +  \int | \partial_x u_{-} |^2 \right) + 2C_3 \|u\|
\]
from which the second inequality in the Lemma follows.
\qed

We now begin to assemble the pieces to prove Proposition \ref{prop:harmonic}.
From the definition of $\Cell$, 
\begin{eqnarray*}
\frac{1}{8} \| \partial_x u \|^2 + \frac{\beta_0}{2\nu} \| \Cell u \|^2  
&\ge & C_1 \| \partial_x u \|^2 + \frac{C_2l^2\beta_0e^{-2\nu t}}{\nu} \| \cos x\ u \|^2,
\end{eqnarray*}

We now use Lemma \ref{lem:plus_minus} to bound this expression from below by
\[
\left(  \frac{2 C_1}{3}  \| \partial_x u_{+} \|^2 + \frac{C_2l^2\beta_0e^{-2\nu t}}{ \nu} \| \cos x\ u_{+} \|^2 \right) 
+  \left( \frac{2 C_1}{3}  \| \partial_x u_{-} \|^2 + \frac{C_2l^2\beta_0e^{-2\nu t}}{ \nu} \| \cos x\ u_{-} \|^2 \right)
- C_3 \| u \|^2
\]

\begin{Lemma} \label{lem:gallay} There exists $M_4$, such that for any $u_{+}$ with support
in $[-\delta, \pi+ \delta ]$, 
\[
 \frac{2 C_1}{3}  \| \partial_x u_{+} \|^2 + \frac{C_2l^2\beta_0e^{-2\nu t}}{ \nu} \| \cos x\ u_{+} \|^2 \ge
 \frac{M_4 | \ell |\sqrt{\beta_0}e^{-\nu t}}{\sqrt{\nu} } \| u_{+} \|^2\ .
\]
\end{Lemma}
\proof 
There exists $D > 0$ such that for any $-\delta < x < \pi+\delta$, $\cos^2 x \ge D (x-\pi/2)^2$.
Thus
\[
\frac{2 C_1}{3}  \| \partial_x u_{+} \|^2 + \frac{C_2l^2\beta_0e^{-2\nu t}}{ \nu} \| \cos x\ u_{+} \|^2
\ge \int \left( \frac{2 C_1}{3} | \partial_x u_{+} |^2 + \frac{D C_2l^2\beta_0e^{-2\nu t}}{ \nu} (x - \pi/2)^2 | u_{+} |^2 \right).
\]
But, as pointed out in \cite[Appendix A]{GallagherGallayNier09}, this last integral is the quadratic form associated to the quantum mechanical harmonic oscillator. In more standard notation, 
\[
H = -h^2 \partial_x^2 + \omega^2x^2 \qquad \Rightarrow \qquad (Hu, u)_{L^2(\mathbb{R})} 
\geq \hbar\omega(u,u)_{L^2(\mathbb{R})}.
\]
Transferring this to the current setting, we obtain the lemma.
\qed

If we also apply Lemma \ref{lem:gallay} to $u_{-}$, we obtain
\[
\frac{1}{8} \| \partial_x u \|^2 + \frac{\beta_0}{2\nu} \| C u \|^2  \ge \frac{M_4 | \ell |\sqrt{\beta_0}e^{-\nu t}}{\sqrt{\nu}} \left( \| u_{+} \|^2
+ \| u_{-} \|^2 \right) - C_3 \| u \|^2
\]
Now apply Lemma \ref{lem:plus_minus} and find that 
\[
\frac{1}{8} \| \partial_x u \|^2 + \frac{\beta_0}{2\nu} \| \Cell u \|^2  \ge \frac{3 M_4| \ell |\sqrt{\beta_0}e^{-\nu t}}{4\sqrt{ \nu}} \| u \|^2 - C_4 \| u \|^2
\]
from which Proposition \ref{prop:harmonic} follows if $\nu$ is sufficiently small.
\end{Proof}

We conclude this section by briefly discussing how we might extend the theoretical results of this
 section to the full linearized operator (\ref{E:linear}).  Villani's hypocoercivity method applies to operators of
 form
 \begin{equation}
{\mathcal L} = A^* A + B
\end{equation}
where $B$ is anti-symmetric.  If we consider the linear operator  (\ref{E:linear}) restricted to functions
with fixed Fourier coefficient $l$, as we did for the approximation operator earlier in this section
it can be written as
\begin{equation}\label{E:Lell}
\cLell \hwell = \nu \Delta_{\ell} \hwell -i \ell a e^{-\nu t} (\sin x) (1+ \Delta_{\ell}^{-1}  ) \hwell\ ,
\end{equation}
where $\Delta_{\ell} = (\partial_x^2 - \ell^2)$.  We can write $\nu  \Delta_{\ell}$ as $A^* A$ (as we did
for the approximation linearization studied here), but the remainder of the operator is no longer anti-symmetric.

However, note that on the space ${\mathcal M}$ defined in Section \ref{S:anomolous_modes}, (i.e. the complement of the
``slow'' modes), the operator $(1+\Delta_{\ell}^{-1})$ is a symmetric, strictly positive operator, and hence we can
define its square root, $\sqrt{1+\Delta_{\ell}^{-1} }$.  Then, if we make the change of variables
$\huell = \sqrt{1+\Delta_{\ell}^{-1} } \hwell$, this transforms the operator $\cLell \hwell$ to
\begin{equation}
\tilde{\cLell} \huell = \nu \Delta \huell -i \ell a e^{-\nu t} \left[ \sqrt{1+\Delta_{\ell}^{-1} } \sin x \sqrt{1+\Delta_{\ell}^{-1} } \right] \huell
\end{equation}
The operator $\tilde{\Bell} = -i \ell a e^{-\nu t} \left[ \sqrt{1+\Delta_{\ell}^{-1} } \sin x \sqrt{1+\Delta_{\ell}^{-1} } \right] $
is anti-symmetric, so we can attempt to apply Villani's theory to this new form of the operator.
We note that in passing
having written the operator in this form (i.e. the sum of a dissipative operator, plus an anti-symmetric operator),
we immediately obtain that the bar states are stable because all eigenvalues of $\tilde{\cLell}$ must have negative
real part.  However, this abstract argument gives only a weak estimate on how large these negative parts are - the
hypocoercivity argument has the potential to give much stronger estimates on the location of these eigenvalues.

If we define
\begin{equation}
\tilde{\Cell} = [\partial_x, \tilde{\Bell}] =  -i \ell a e^{-\nu t} \left[ \sqrt{1+\Delta_{\ell}^{-1} } \cos x \sqrt{1+\Delta_{\ell}^{-1} } \right]
\end{equation}
most of the preceding argument goes through without change.  However the fact that 
$[\tilde{\Bell}, \tilde{\Cell}] \ne 0$ causes problems at several points.  We believe that this problem can be overcome
by altering the functional $\Phi^{\ell}(t)$, and plan to address this point in future research.  We note that similar ideas
have recently been used by Deng \cite{Deng12} to study the  linearization of the two-dimensional Navier-Stokes
equation about the Oseen vortex solutions. 


\section{Numerics}\label{S:numerics}

The results of the previous section prove that for our approximation to the linearization about the bar states,
solutions decay at a rate of at least ${\mathcal{O}}(e^{- \sqrt{\nu} t})$, much faster than the viscous
decay rate of ${\mathcal{O}}(e^{- \nu t})$.
In this section we present numerical evidence that this rate of approach is actually optimal, not just for
the approximate linearization studied in the previous section, but also for the full
linearization \eqref{E:linear}.  We numerically
compute the eigenvalues of the linearization of the two-dimensional vorticity equation about the bar states,
for a fixed value of the time $t$, and show that as we vary the viscosity parameter $\nu$, the real part of these
eigenvalues are negative and
scale like $ c  \sqrt{\nu} $, so long as we avoid the anomalous subspace.

Recall that as we remarked earlier, if we expand the space $L^2([-\pi,\pi]\times[-\pi,\pi])$, with respect to the basis
$e^{i k x} e^{i l y}$, then the operator ${\mathcal L}(t)$ leaves  invariant the subspaces with $ \ell$ fixed.  Therefore we
can study the eigenvalue asymptotics, by fixing $\ell $ and studying the spectrum of the operator
${\mathcal{L}}_{\ell}(t)$ defined in \eqref{E:Lell} for different values of $\ell$.  We compute (approximately) the
spectrum of ${\mathcal{L}}_{\ell}(t)$ by writing out its matrix representation with respect to $\{ e^{i kx } \}_{k=-N}^N$ and
then using Mathematica to compute the eigenvalues of the resulting matrix.  We have examined matrices with
different values of $N$, up to $N=400$ and we find consistent results for all values of $N$ provided $N$ is large
enough (bigger than about $30$-$40$) to insure that we have good asymptotics.

We present two examples of our numerical computations below, both with $N=400$ and $\ell=2$.  Also, we fix the
time $t$ so that $a e^{-\nu t} = 1$.  

The first figure plots the logarithm of the absolute value of the real part of the least negative eigenvalue, 
as a function of the logarithm of $\nu$, for values of $\nu =$ $0.005$, $0.002$, $0.001$, $0.0005$, $0.00025$ and $0.0001$.
Superimposed on these six data points is a line of slope $2$, which corresponds to the real part of the eigenvalue
scaling as $\sqrt{\nu}$.  The fit is quite good, especially for the smaller values of $\nu$.

\begin{figure}[htp]
\centering
\includegraphics[width=6in,height=3in]{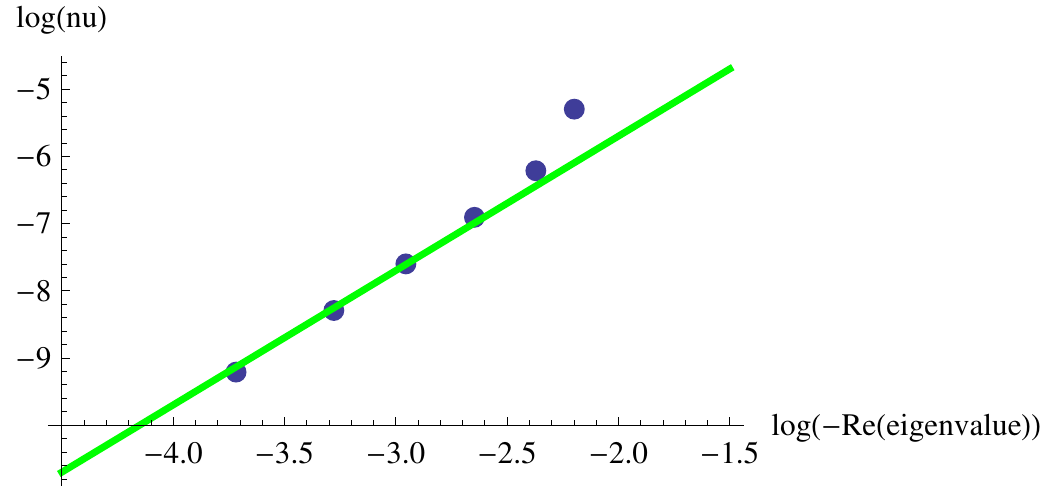}
\caption{The computed value of the real part of the eigenvalue of the linearization about
a bar state, compared with a line which represents a scaling proportional to $\sqrt{\nu}$.}
\label{fig:single_eigenvalue}
\end{figure}

\begin{figure}[htp]
\centering
\includegraphics[width=6in,height=3in]{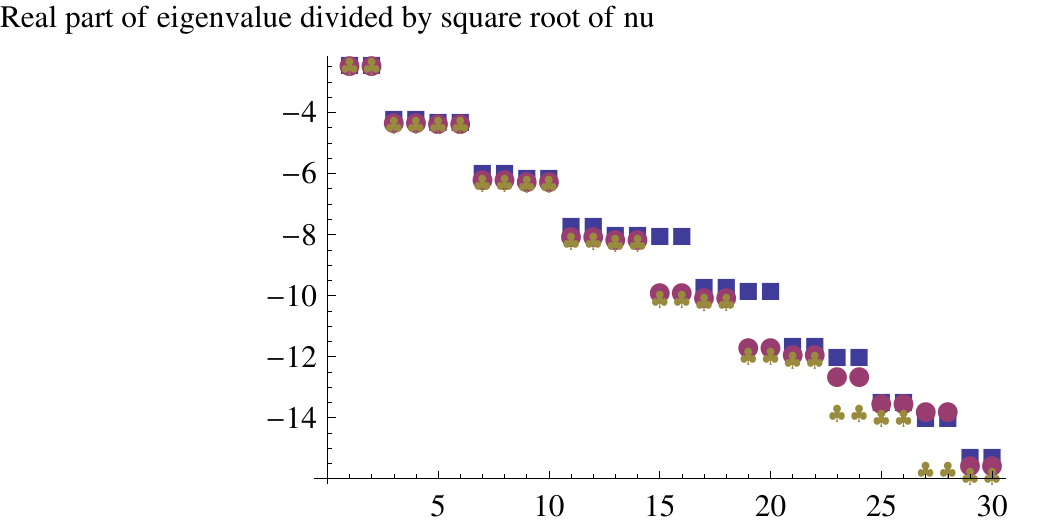}
\caption{The scaling behavior of the first thirty eigenvalues of the linearization about the bar states, with $l=2$ and
a $801 \times 801$ matrix approximation to the $x$-dependent part of the operator.  In this figure, solid squares correspond
to the eigenvalues with $\nu=0.00025$, solid circles correspond to $\nu = 0.0001$ and clubs correspond to $\nu = 0.00005$. }
\label{fig:scaling}
\end{figure}

For our second piece of numerical evidence in support of the optimality of the $\sqrt{\nu}$ scaling of
the attractivity of the bar states, we first compute the eigenvalues of linearization about the bar states, and
then order the resulting eigenvalues according to decreasing real part.  We consider the first $30$ eigenvalues
in this ordering and plot the real part of each eigenvalue, divided by $\sqrt{\nu}$.  We repeat this procedure
for three different values of $\nu$, $\nu=0.00025$, $\nu=0.0001$ and $\nu=0.00005$.  If the eigenvalues were
of the form $C_k \sqrt{\nu}$, $j=k, \dots , 30$, then these three data sets would fall on top of each other.
As Figure \ref{fig:scaling} shows, this is nearly true, especially for the first few eigenvalues.  As one moves
to eigenvalues further down the list, the smaller values of $\nu$ overlap more than the larger values.  This
is to be expected - for larger values of $k$, the diagonal part of the linearization $-\nu (l^2+k^2)$ will dominate
the behavior of the eigenvalue until $\nu$ becomes quite small  and the interaction between the diagonal
and off-diagonal parts captured by the hypocoercivity method makes itself fully felt.

\section{Discussion} \label{S:discussion}

We have analyzed the linearization of the two-dimensional vorticity equation (\ref{eq:2Dvort}) about the bar states.
The bar states are ``quasi-stationary''  states since for small values of the viscosity, $\nu$, they decay 
on the very slow scale $\mathcal{O}(e^{-\nu t})$.
 We have proposed a new, dynamical systems based, explanation for the metastable behavior these states
exhibit in numerical simulations.  Our explanation begins by  showing that a natural approximation to the 
linearization of the Navier-Stokes equations about the bar states exhibits decay on a much faster time scale than the 
viscous decay rate and also presenting numerical computations of the eigenvalues of the full linearization of the
equation which are consistent with the decay rate we obtain rigorously for our approximation.

Ideally, one would like to extend this linear result to the nonlinear equation (\ref{E:lin_b}).  We are currently working on
this extension, first of all by showing that our theoretical results on the approximate linearization can be extended to
the full, linearized equations.  The direct extension of the linear results to the full, nonlinear equation is complicated by
the fact that the linearized operator is time dependent, so that we can't use ordinary spectral projections to project
the nonlinear equation onto rapidly decaying and slowly decaying modes (the ``anomalous'' modes of the linear
equation, described in Section \ref{S:anomolous_modes})  as is done in the stable manifold theorem, for instance.

Another interesting direction to pursue is to analyze the stability of the dipole states, given explicitly by
\begin{equation}\label{E:dipole}
\omega^{\mathrm{d}, m}(x,y, t) = e^{-\nu m^2 t}[ \cos mx + \cos my], \quad \tilde{\omega}^{\mathrm{d}, m}(x,y, t) = e^{-\nu m^2 t}[ \sin mx + \sin my]
\end{equation}
with associated velocity fields
\[
\uu^{d,m}(x,y,t) =  \frac{1}{m} e^{-\nu m^2 t} \begin{pmatrix} -\sin my \\ \sin mx \end{pmatrix}, \qquad \tilde{\uu}^{d,m}(x,y,t) = - \frac{1}{m} e^{-\nu m^2 t} \begin{pmatrix} -\cos my \\ \cos mx \end{pmatrix}.
\]
These states are also maximum entropy solutions of the Euler equations, and hence likely candidates for relevant quasi-stationary states for the 2D vorticity equation. In fact, these states seem to be more important than the bar states, at least for a square torus \cite{BouchetSimonnet09}, as solutions with a larger class of initial data converge to them. See \cite[Figures 2a and 5]{YinMontgomeryClercx03}. However, we have not yet found an invariant subspace of solutions that converge rapidly to these states. 

If one linearizes about the dipole solution $a\omega^{\mathrm{d}, 1}$, one finds
\begin{eqnarray}
\partial_t w &=& \mathcal{L}^{d}(t) w -  \vv \cdot \nabla w \label{E:lin_tg} \\
\mathcal{L}^{d}(t) w &=& \nu\Delta w - ae^{-\nu t} [ \sin x \partial_y (1 + \Delta^{-1}) - \sin y \partial_x(1+\Delta^{-1})] w. \nonumber
\end{eqnarray}
In Fourier space, this operator becomes
\begin{eqnarray}
\hat{\mathcal{L}}^{tg}(t) \hat{w}(k,l) &=& - \nu(k^2+l^2) \hat{w}(k,l) -  \frac{l}{2}ae^{-\nu t} \Bigg[ \left( 1 - \frac{1}{(k-1)^2+l^2}\right) \hat{w}(k-1,l) \nonumber \\
&\qquad& \qquad - \left( 1 - \frac{1}{(k+1)^2+l^2}\right) \hat{w}(k+1,l) \Bigg]  \nonumber \\
 &&\hspace*{-.3in}+ \frac{k}{2}ae^{-\nu t} \Bigg[ \left( 1 - \frac{1}{k^2+(l-1)^2}\right) \hat{w}(k,l-1)  - \left( 1 - \frac{1}{k^2+(l+1)^2}\right) \hat{w}(k,l+1) \Bigg] 
\label{E:tg_lin}
\end{eqnarray}

Although $e^{\pm \rmi x}$ and $e^{\pm \rmi y}$ are still eigenfunctions, the adjoint eigenfunctions seem to be more difficult to determine. (An explicit one is $\cos x + \cos y$.) A similar change of variables can be used to transform the second part of the linear operator into a skew symmetric one: by defining $w = \sqrt{1 + \Delta^{-1}} u$, the linear part of equation (\ref{E:lin_tg}) becomes
\begin{equation} \label{E:tg_as}
u_t = \nu \Delta u + ae^{-\nu t} Q u, \quad Q = \sqrt{1 + \Delta^{-1}}[\sin y \partial_x -\sin x \partial_y] \sqrt{1 + \Delta^{-1}},
\end{equation}
where $Q$ is anti-symmetric.  While we have some preliminary numerical data indicating that
most modes again decay at a rate much faster than the viscous time scale, we have not yet proven a result similar to Theorem \ref{Thm:bar_conv} for this operator.


\section{Acknowledgements} The authors wish to thank J.~P.~Eckmann for fruitful discussions, which lead to the discovery of the transformation creating a skew-symmetric operator.  The authors also thank Th.~Gallay for a very helpful explanation of the results and methods used in 
\cite{GallagherGallayNier09}.  M.B. wishes to thank Edgar Knobloch, who pointed out the potential connection between this result and Taylor dispersion. C.E.W wishes to thank G.~van Baalen for helpful discussions of the role of bar states in two-dimensional flows. M. B. was partially funded through NSF grant DMS-1007450. C.E.W. was partially funded by NSF grant DMS-0908093.

\bibliographystyle{plain}
\bibliography{ref_bar}

 
\end{document}